\documentclass[12pt,a4paper]{amsart}
\usepackage{amssymb,amsmath,amsthm}
\usepackage{graphicx}
\usepackage{hyperref}
 \usepackage[color=green!40]{todonotes}

\newcommand\cC{{\mathcal C}}

\makeatletter
\newtheorem*{rep@theorem}{\rep@title}
\newcommand{\newreptheorem}[2]{%
\newenvironment{rep#1}[1]{%
 \def\rep@title{#2 \ref{##1}}%
 \begin{rep@theorem}}%
 {\end{rep@theorem}}}
\makeatother

\newtheorem{thm}{Theorem}
\newtheorem{theorem}{Theorem}

\newtheorem{lemma}[thm]{Lemma}

\newtheorem{conjecture}[thm]{Conjecture}
\newtheorem{proposition}[thm]{Proposition}

\newtheorem{defn}[thm]{Definition}

\newtheorem{corollary}[thm]{Corollary}

\newtheorem{example}[thm]{Example}

\newtheorem{clm}[thm]{Claim}

\newcommand\cref[1]{Corollary~\ref{cor:#1}}

\title{Counting Connected Partitions of Graphs}

\author{Yair Caro}
\address{Department of Mathematics, University of Haifa-Oranim, Israel}
\email{yacaro@kvgeva.org.il}

\author{Bal\'azs Patk\'os}
\address{Alfr\'ed R\'enyi Institute of Mathematics}
\email{patkos@renyi.hu}
\thanks{Research of Patk\'os, Tuza, and Vizer is partially supported by NKFIH grants SNN 129364 and FK 132060.}

\author{Zsolt Tuza}
\address{Alfr\'ed R\'enyi Institute of Mathematics and University of Pannonia}
\email{tuza.zsolt@mik.uni-pannon.hu}

\author{M\'at\'e Vizer}
\address{Alfr\'ed R\'enyi Institute of Mathematics}
\email{vizermate@gmail.com}

\date{}
\date{}

\begin{document}
\begin{abstract}
    
    Motivated by the theorem of Gy\H ori and Lov\'asz, we consider the following problem. For a connected graph $G$ on $n$ vertices and $m$ edges determine the number $P(G,k)$ of unordered solutions of positive integers $\sum_{i=1}^k m_i = m$ such that every $m_i$ is realized by a connected subgraph $H_i$ of $G$ with $m_i$ edges such that $\cup_{i=1}^kE(H_i)=E(G)$. We also consider the vertex-partition analogue.
    
     We prove various lower bounds on $P(G,k)$ as a function of the number $n$ of vertices in $G$, as a function of the average degree $d$ of $G$, and also as the size $\mathrm{CMC}_r(G)$ of $r$-partite connected maximum cuts of $G$. Those three lower bounds are tight up to a multiplicative constant. 
     
     We also prove that  the number $\pi(G,k)$ of unordered $k$-tuples with $\sum_{i=1}^kn_i=n$, that are realizable by vertex partitions into  $k$ connected parts of respective sizes $n_1,n_2,\dots,n_k$, is $\Omega(d^{k-1})$. 
\end{abstract}

\maketitle

\section{Introduction}

Partitions and decompositions of graphs are fundamental subjects in graph theory with a long history.
One of the most well-known partition theorems is the following result from the 1970s, which was independently proved by Lov\'asz \cite{L} and Gy\H ori \cite{Gy} (see also \cite{HT}); it serves as the basic motivation of our paper.

\begin{theorem}[Gy\H ori \cite{Gy}, Lov\'asz \cite{L}]\label{GyLvertex}
Let\/ $k \ge 2$ be an integer. If\/ $G$ is a\/ $k$-connected graph on\/ $n$ vertices, $\sum_{i=1}^kn_i=n$ with integers\/ $n_i\ge 1$, and\/ $v_1, \ldots, v_k \in V(G)$ are distinct vertices, then there exist\/ $V_1,V_2,\dots,V_k $ with\/ $|V_i|=n_i$, $v_i \in V_i$, $V(G)=\cup_{i=1}^kV_i$ such that the induced subgraph\/ $G[V_i]$ is connected for all\/ $i=1,2,\dots,k$.
\end{theorem}

The proofs of the full strength of Theorem \ref{GyLvertex} are non-constructive. For recent developments on algorithmic aspects of the problem and some related results see \cite{BL,Betal,CFINZ22} and the references therein.

Note that if $k\ge2$ and $G$ is $k$-edge-connected, then $L(G)$, the line graph of $G$ is $k$-connected (see  \cite{CS69}). Hence we can state the following version of Theorem \ref{GyLvertex} for $k$-edge-connected graphs. 

\begin{theorem}\label{GyLedge}
Let\/ $k \ge 2$ be an integer. If\/ $G$ is a\/ $k$-edge-connected graph with\/ $m$ edges, $\sum_{i=1}^km_i=m$ with integers\/ $m_i\ge 1$, and\/ $e_1, \ldots, e_k \in E(G)$ are distinct edges, then there exist\/ $E_1,E_2,\dots,E_k $ with\/ $|E_i|=m_i$, $e_i \in E_i$, $E(G)=\cup_{i=1}^kE_i$ such that the subgraph induced by the edges of\/ $E_i$ is connected for all\/ $i=1,2,\dots,k$.

\end{theorem}

An immediate consequence of Theorem \ref{GyLvertex} and Theorem \ref{GyLedge} is that the number of different unordered $k$-tuples that we can achieve as sizes of the parts in a partition of a $k$-edge-connected (resp., $k$-vertex-connected) graph into $k$ edge-disjoint (resp., $k$ vertex-disjoint) connected subgraphs is exactly the number of different unordered $k$-tuples of positive integers whose sum is $m=|E(G)|$ (resp., $n=|V(G)|$).

This is in fact the number-theoretic \textit{partition function} $\pi(n,k)$ that counts the number of representations of an integer $n$ as a sum of $k$ positive integers.
Its asymptotics was determined by Erd\H os and Lehner \cite{EL} who showed that $\pi(n,k)=(1-o(1))\frac{\binom{n-1}{k-1}}{k!}$ holds as long as $k=o(n^{\frac{1}{2}})$.

Our main goal in this article is to estimate for $k \ge 2$ the number of different unordered $k$-tuples that represent the sizes (number of edges) of the parts in a partition of the edge set $E(G)$ of a connected graph into $k$ edge-disjoint connected subgraphs. We also deal with the vertex-partition analogue. We shall use the following notation.

For a graph $G$ we use $n=n(G)$ for the number of vertices and $m=e(G)$ for the number of edges. 

\begin{defn}
For a graph\/ $G$, a connected\/ $(m_1,m_2,\dots,m_k)$-edge partition with\/ $\sum_{i=1}^km_i=e(G)$ is a partition\/ $E(G)=E_1\cup E_2 \cup \dots \cup E_k$ with\/ $|E_i|=m_i$ for all\/ $i=1,2,\dots,k$ such that every\/ $E_i$ forms a connected graph on the vertex set\/ $\cup_{e\in E_i}e$. For a connected graph\/ $G$ and an integer\/ $k\ge 2$ let\/ $$P(G,k):= |\{(m_1,m_2,\dots,m_k): m_1\ge m_2 \ge \dots \ge m_k \ge 1, \ $$ $$\qquad\qquad\qquad\qquad \exists \ \textrm{connected} \ (m_1,m_2,\dots,m_k)\textrm{-edge partition of} \ G\}|.$$ 
\end{defn}

\begin{defn}
For a graph\/ $G$, a connected\/ $(n_1,n_2,\dots,n_k)$-vertex partition with\/ $\sum_{i=1}^kn_i=n(G)$ is a partition\/ $V(G)=V_1\cup V_2 \cup \dots \cup V_k$ with\/ $|V_i|=n_i$ for all\/ $i=1,2,\dots,k$ such that every\/ $V_i$ induces a connected subgraph. For a connected graph\/ $G$ and an integer\/ $k\ge 2$ let\/ $$\pi(G,k):= |\{(n_1,n_2,\dots,n_k): n_1\ge n_2 \ge \dots \ge n_k \ge 1, \ $$ $$\qquad\qquad\qquad\qquad \exists \ \textrm{connected} \ (n_1,n_2,\dots,n_k)\textrm{-vertex partition of} \ G\}|.$$ 
\end{defn}

\

So the basic problems that we address in this article are: 
\begin{itemize}
    \item Find lower bounds for $P(G,k)$.
    \item Find lower bounds for $\pi(G,k)$.
\end{itemize}

As we shall see later, such estimates can be obtained depending upon the average degree $d=d(G) := 2m/n$.  

\medskip

For some results concerning decompositions of graphs into isomorphic subgraphs, which are related to our edge-partition problem, see e.g.\ \cite{DF85,JRP85}.

We use several well-known results concerning the existence of $k$-connected subgraphs and the existence of $k$ edge-disjoint spanning trees, which we state now for later use.  

\begin{theorem}[Mader \cite{M72}]\label{mader} Let\/ $\ell \ge 1$. If\/ $d(G) \ge 4\ell$, then $G$ contains an\/ $(\ell+1)$-connected subgraph\/ $H$ with\/ $d(H) > d(G) - 2\ell$.

\end{theorem}

We shall comment on recent improvements of Mader's theorem in the concluding section, but for our purpose Theorem \ref{mader} suffices.

\begin{theorem}[Tutte \cite{T61}, Nash-Williams \cite{NW61}]\label{tutte} Every finite\/ $2k$-edge-connected graph has\/ $k$ edge-disjoint spanning trees.

\end{theorem}

Lastly, we shall apply \emph{connected maximum cut} (CMC in short), which is a variant of the Max-Cut problem. More precisely let $V(G) = A \cup B$ be a partition of V(G), and let $e(A, B)$ denote the number of crossing edges between the parts (i.e.\ edges with one end in $A$ and the other end in $B$). Then let us define $$\mathrm{CMC}(G) := \max \{  e(A,B) :  V(G) = A \cup  B, \ \textrm{both} \ A \ \textrm{and} \ B \ \textrm{induce connected subgraphs} \}.$$  

As far as we are aware, unlike the much work that was done on estimating Max-Cut in terms of $e(G)$ for example (see e.g.\ \cite{LNS16,H14} or the extensive survey \cite{PT95} and references therein), $CMC$ seems to have been the subject of systematic research only in the last decade, and only from the algorithmic point of view, seeking good approximations for $\mathrm{CMC}(G)$; see \cite{HKMPS15,HKMPS20}.

We shall also generalize the result on $\mathrm{CMC}(G)$ to the more general setting concerning the maximum number of edges in a spanning $r$-partite subgraph whose color classes (i.e., parts) are connected, noting that $\mathrm{CMC}(G)$ is just the case $r=2$. We shall denote this parameter by $\mathrm{CMC}_r(G)$. 

Before presenting our results, let us mention how $\mathrm{CMC}(G)$ and $P(G,2)$ are related. Suppose $A,B$ is a maximum connected cut of $G$ i.e.\ $V(G)=A\cup B$ and $G[A],G[B]$ are both connected and $e(A,B)=\mathrm{CMC}(G)$. Then if $e_1,e_2,\dots,e_m$ is an enumeration of all edges from $A$ to $B$, then $E^1_j=E(G[A])\cup \{e_1,e_2\dots,e_j\}$ and $E^2_j=E(G[B])\cup \{e_{j+1},e_{j+2},\dots,e_m\}$ form a connected edge-partition of $G$ such that $(|E^1_j|,|E^2_j)\neq (|E^1_h|,|E^2_h|)$ if $j\neq h$. Therefore $P(G,2)\ge \frac{\mathrm{CMC}(G)}{2}$ holds for any graph $G$.

\subsection{Our results}

In our first theorem, we give a lower bound on $P(G,k)$, if $G$ is a connected graph. 

\begin{thm}\label{kmain}
For every\/ $k\ge 2$ there exists a constant\/ $c_k>0$ such that for any connected graph\/ $G$ on\/ $n$ vertices we have\/ $P(G,k)\ge c_k\log^{k-1}n$. Moreover, if\/ $n=4\cdot 3^{\ell}+\sum_{i=0}^{\ell}3^i$ for some\/ $\ell\ge 1$, then for any tree\/ $T$ on\/ $n$ vertices we have\/ $P(T,2)\ge P(T(\ell),2)$, where\/ $T(\ell)$ denotes the tree on $4\cdot 3^{\ell}+\sum_{i=0}^{\ell}3^i$ vertices that we obtain from the complete ternary tree of height \/$\ell$ such that we identify every leaf with the middle vertex of a path on 5 vertices. 
\end{thm}

We obtain the following lower bounds on $P(G,k)$ and $\mathrm{CMC}(G)$ with respect to average degree $d=d(G)$. 

\begin{thm}\label{lowerpg}
For any connected graph\/ $G$ with average degree\/ $d$, $P(G,k)\ge  \Omega_k(d^{2(k-1)})$.
\end{thm}

\begin{thm}\label{cmcg}
For any connected graph\/ $G$ with average degree\/ $d$, $\mathrm{CMC}_r(G)\ge \Omega_r(d^2)$.
\end{thm}

We shall show in the proof of Theorem \ref{cmcg} that the case $r=2$ supplies another proof for the case $k=2$ of Theorem \ref{lowerpg}. The next theorem shows that Theorems \ref{kmain}, \ref{lowerpg}, and \ref{cmcg} can only be improved by a constant multiplicative factor. To see the tightness of Theorem \ref{cmcg}, one uses the above inequality $P(G,2)\ge \frac{\mathrm{CMC}(G)}{2}$.

\begin{thm}\label{constr}
Let us fix $k\ge 2$ and for $a\ge 1$, we write $n_a=\sum_{j=0}^a2^j$. For any function\/ $2-\frac{2}{n}\le d(a)\le n_a-1$, there exists a sequence of connected graphs\/ $G_a$ on\/ $n_a$ vertices with average degree\/ $d_a$ such that\/  $P(G_a,k)=\Theta_k((\log n_a +d^2(a))^{k-1})$ and if $d(a)$ tends to infinity, then $\lim_{a\rightarrow \infty}\frac{d_a}{d(a)}=1$, otherwise $\frac{d_a}{d(a)}=\Theta(1)$. 
\end{thm}

Putting together Theorems \ref{kmain}, \ref{lowerpg} and \ref{constr}, we obtain the following corollary.

\begin{corollary}
The minimum of\/ $P(G,2)$ over all\/ $n$-vertex connected graphs with average degree\/ $d$ is\/ $\Theta(\log n+ d^2)$.

\end{corollary}

Finally, we have the following theorem about $\pi(G,k)$.

\begin{thm}\label{thmonpi}
Let\/ $G$ be a connected graph on\/ $n$ vertices with average degree\/ $d$. Then for\/ any fixed $k \ge 2$ we have\/ $\pi(G,k) =\Omega_k(d^{k-1})$.
\end{thm}

\section{Proofs}

We start with two auxiliary statements and a definition. Let $t(n)$ be the sequence defined by $t(1)=0$, $t(2)=1$, $t(3)=2$, and $t(n)=\min\{d+t(\lceil\frac{n-1}{d}\rceil):d\ge 1\}$ for all $n\ge 4$.

\begin{proposition}\label{t}
The sequence\/ $t(n)$ 
\begin{itemize}
    \item 
    is monotone increasing,
    \item 
    if $n \ge 11$, then $t(n)=3+t(\lceil \frac{n-1}{3}\rceil)$,
    \item
    satisfies\/ $t(n)\ge 3\log_3n-C$ for some absolute constant\/ $C$,
    \item
    satisfies\/ $t(n)=P(T(\ell),2)+2$ for all\/ $n=4\cdot 3^\ell+\sum_{i=0}^{\ell}3^i$ with $\ell\ge 1$.
\end{itemize}
\end{proposition}

\begin{proof}
Monotonicity follows by induction easily as it is true for $n=1,2,3$. 

Observe that by the monotonicity of $t$, the set $t^{-1}(h)=\{n: t(n)=h\}$ is an interval for any $h$.

Then to obtain the  equality $t(n)=3+t(\lceil \frac{n-1}{3}\rceil)$, it is enough to check that for any $d$ and $h$, if the smallest element of $t^{-1}(h)$
is $n$, then $\lceil \frac{n-1}{d}\rceil$ belongs to $t^{-1}(h')$ for some $h'\ge h-d$. We claim that this is true if $d\ge 4$. Indeed, for any $d\ge 4$ there exist $d_1,d_2$ such that $d_1+d_2\le d$ but $d_1\cdot d_2\ge d$. For $d=5,6$, the pair $d_1=2,d_2=3$ works, while for $d\ge 7$ and $d=4$ the pair $d_1=d_2=\lceil \sqrt{d}\rceil$ works. Then 
\[
t(n)\le t(\lceil\frac{n-1}{d_1}\rceil)+d_1\le t(\lceil \frac{\lceil\frac{n-1}{d_1}\rceil-1}{d_2}\rceil)+d_1+d_2\le t(\lceil\frac{n-1}{d}\rceil)+d.
\]
The statement for $d=1$ is simply monotonicity. We are left to show the case $d=2$. We do this by determining the intervals $t^{-1}(h)$ for all $h$. The statements $t^{-1}(0)=1, t^{-1}(1)=2, t^{-1}(2)=3, t^{-1}(3)=\{4,5\}, t^{-1}(4)=\{6,7\}, t^{-1}(5)=[8,11], t^{-1}(6)=[12,16],t^{-1}(7)=[17,23]$ can be checked by hand. Then we claim
\begin{itemize}
    \item 
    if $h = 3k-1$, then $t^{-1}(h)=[\frac{3^{k} - 1}{2} + 3^{k-1} + 3^{k-3}+1, \frac{3^{k+1} - 1}{2}  - 3^{k-1} + 3^{k-2}]$,
    \item 
    if $h = 3k$, then $t^{-1}(h)=[\frac{3^{k+1} - 1}{2}  - 3^{k-1} + 3^{k-2}+1, \frac{3^{k+1} - 1}{2}  + 3^{k-1}]$,
    \item
    if $h = 3k+1$, then $t^{-1}(h)=[\frac{3^{k+1} - 1}{2}  + 3^{k-1}+1, \frac{3^{k+1} - 1}{2} + 3^k + 3^{k-2}]$.
\end{itemize}
This again can be verified for $h=8,9,10$ by hand.  For the general statement for $d=2$, we need to show that for any $h$ if the smallest element of $t^{-1}(h)$
is $n$, then $\lceil \frac{n-1}{2}\rceil$ belongs to $t^{-1}(h')$ for some $h'\ge h-2$. That is $\lceil \frac{n-1}{2}\rceil\ge n'$, where $n'$ is the smallest element of $t^{-1}(h-2)$. We consider the cases $h=3k-1$, $h=3k$, and $h=3k+1$.

\begin{itemize}
    \item 
    Suppose first $h=3k-1$. Then  we need
    \[
    \lceil \frac{\frac{3^{k} - 1}{2} + 3^{k-1} + 3^{k-3}}{2}\rceil=\lceil (11+\frac{3}{4})\cdot 3^{k-3}-\frac{1}{4}\rceil \ge (11+\frac{1}{2})\cdot3^{k-3}+\frac{1}{2}=\frac{3^{k} - 1}{2}  - 3^{k-2} + 3^{k-3}+1
    \]
    if $k\ge 4$ that is $h\ge 11$.


    \item 
     Suppose next $h=3k$. Then  we have $h-2=3k-2=3(k-1)+1$ and so we need
      \[
    \lceil \frac{\frac{3^{k+1} - 1}{2}  - 3^{k-1} + 3^{k-2}}{2}\rceil=\lceil \frac{23}{4}\cdot 3^{k-2} -\frac{1}{4}\rceil \ge \frac{11}{2}\cdot3^{k-2}+\frac{1}{2}=\frac{3^{k} - 1}{2}  + 3^{k-2}+1
    \]
    if $k\ge 4$ that is $h\ge 12$.


    \item 
    Suppose finally $h=3k+1$. Then  we have $h-2=3k-1$ and so we need
      \[
    \lceil \frac{\frac{3^{k+1} - 1}{2} + 3^{k-1} }{2}\rceil=\lceil \frac{99}{4}\cdot 3^{k-3}-\frac{1}{4}\rceil \ge \frac{47}{2}\cdot 3^{k-3}+\frac{1}{2}=\frac{3^{k} - 1}{2} + 3^{k-1} + 3^{k-3}+1
    \]
    if $k\ge 4$ that is $h\ge 13$.


\end{itemize}

Now $t(n)\ge 3\log_3n-C$ follows from monotonicity and from the fact that the above argument implies that for $n=10\cdot \sum_{i=0}^\ell3^i$ we have $t(n)=t(10)+3\ell$.

To see the last statement, one only needs the trivial observations $P(T(\ell+1),2)=3+P(T(\ell),2)$ and $P(T(1),2)=4=t(16)-2$.
\end{proof}

\begin{lemma}\label{nested}
Let\/ $T$ be a tree on\/ $n$ vertices and\/ $v\in V(T)$. Then there exists a sequence\/ $(A_1,B_1), (A_2,B_2),\dots,(A_m,B_m)$ such that 
\begin{enumerate}
    \item 
    $A_i,B_i \subseteq V(T)$, $A_i\cap B_i= \{v_i\}$, $A_i\cup B_i=V(T)$, $v_1=v$, and for any\/ $i<j$, $v_i\in B_j$,
    \item
    $T[A_i], T[B_i]$ are connected for all\/ $i=1,2,\dots,m$,
    \item
    $A_1\supsetneq A_2 \supsetneq \dots \supsetneq A_m$ and\/  $B_1\subsetneq B_2 \subsetneq \dots \subsetneq B_m$,
    \item
    $m\ge t(n)+1$.
\end{enumerate}
\end{lemma}

\begin{proof}
We apply induction on $n$. For $n=1,2,3$, the statement trivially holds.  Let the components of $T\setminus \{v\}$ be $D_1,D_2,\dots,D_d$ with $|D_d|\ge (n-1)/d$, and let $v'$ be the neighbor of $v$ in $D_d$. Then for $i=1,2,\dots,d$ we set $A_i:=V(T)\setminus \cup_{j=1}^{i-1}D_j$ and $B_i=\{v\}\cup \bigcup_{j=1}^{i-1}D_j$ (so $A_1=V(T)$ and $B_1=\{v\}$). Clearly, (2) and (3) are satisfied so far. We apply induction to $T[D_d]$ and $v'$, to obtain a sequence $(A'_{j},B'_{j},u_{j})$ for $j=1,2,\dots,m'$ with $m' \ge t(\lceil \frac{n-1}{d}\rceil)+1$. As all $B'_{j}$ contain $v'$, $T[B'_{j}\cup (V(T)\setminus D_d)]$ is connected. Therefore, we can set $A_{d+j}=A'_{j}$, $B_{d+j}=B'_j\cup (V(T)\setminus D_d)$, $v_{d+j}=u_j$, to
get a sequence satisfying (1), (2), and (3). Its length is at least $d+t(\lceil \frac{n-1}{d}\rceil)+1\ge t(n)+1$. We used the definition of $t(n)$ to obtain the last inequality.
\end{proof}

\begin{proof}[\textbf{Proof of Theorem \ref{kmain}}]
We apply induction on $k$. We start with the case $k=2$. Let $G$ be a connected graph on $n$ vertices and let $T$ be a spanning tree of $T$. Applying Lemma \ref{nested}, we obtain a sequence  $(A_i, B_i)$ of length at least $3\log_3n-C$ satisfying properties (1), (2), (3), and (4). Let us define the edge-partitions $(E^i_1,E^i_2)$ as $$E^i_1:=E(G[A_i]) \cup E(G[A_i,B_i]) \hskip 1truecm \text{and} \hskip 1truecm E^i_2=E(G[B_i]).$$
By definition, $E^i_1\cup E^i_2=E(G)$ for all $i$ and $|E^i_2|< |E^{i+1}_2|$. The subgraphs  induced by $E^i_1$ and $E^i_2$ are connected as $E^i_2=E(G[B_i])$ lives on $B_i$ and contains $E(T[B_i])$ which is connected by Lemma \ref{nested}. Also, $E^i_1$ is connected as it contains $E(G[A_i])\supseteq E(T[A_i])$, which is connected on $A_i$ by Lemma \ref{nested}, and all other edges of $E^i_1$ have one endpoint in $A_i$. So we obtained at least $3\log_3n-C$ connected edge-partitions of $G$ such that the ordered pairs $(|E^i_1|,|E^i_2|)$ are all distinct. Thus we have $P(G,2)\ge \frac{3\log_3n-C}{2}$.

Suppose next that the statement holds for $k$. Let $G$ be a connected graph on $n$ vertices and $T$ a spanning tree of $G$. We plan to apply Lemma \ref{nested} to $T$ to obtain a sequence $(A_i, B_i)$. It is known that in every $T$ on $n$ vertices, there exists a vertex $v$ such that all components of $T-v$ have size at most $n/2$. If $T- v$ contains two components $C_1,C_2$ with $|C_1|\le |C_2|$, then we let $A_1=C_1\cup \{v\},B_1=C_2\cup \{v\}$. If $T-v$ has at least three components, then there is a partition of the components into two parts $\cC_1$, $\cC_2$ such that $n/3 \le \sum_{C\in \cC_1}C,\sum_{C\in \cC_2}C\le 2n/3$, so one of the sums is between $n/3$ and $n/2$. Let it be $\cC_1$ and define $A_1=\{v\}\cup \bigcup_{C\in \cC_1}C,B_1=\{v\}\cup \bigcup_{C\in \cC_2}C$. In both cases, we achieved $A_1$ with $n/3\le |A_i|\le n/2$. We consider the sequence applied to $T[A_1]$ so its length is at least $3\log_3(n/3)-C=3\log_3n-C-3\ge c\log n$ for some $c$. Also, the corresponding $B_i$ all have size at least $n/2$ as $B_1$ has at least that size. Applying the inductive hypothesis, for each $B_i$ we obtain at least $c_k\log^{k-1}(n/2)$ connected edge-partitions of $G[B_i]$. Clearly, $E(G[A_i])\cup E(G[A_i, B_i])$ induces a connected subgraph of $G$ and by the nestedness, their sizes form a decreasing sequence. So, we obtain $(3\log_3(n/3)-C)c_k\log^{k-1}(n/2)$ distinct \textit{ordered} connected edge-partitions. Therefore, we have $$P(G,k+1)\ge \frac{(3\log_3(n/3)-C)c_k\log^{k-1}(n/2)}{(k+1)!}\ge c_{k+1}\log^{k}n.$$

To see the statement about $P(T,2)$, by Proposition \ref{t}, it is enough to prove that for any $n$-vertex tree $T$, we have $P(T,2)\ge t(n)-2$. This can be checked for $n=2,3,4$. Let $n\ge 5$, and let $T$ be an $n$-vertex tree. We consider two cases.

\vskip 0.2truecm

\textsc{Case I.} There exists an edge $uv$ such that $T-uv$ contains two components exactly of the same size.

\vskip 0.15truecm

Let $C_1$ be the component of $T-uv$ that contains $u$, and $C_2$ be the other component. Then we can apply Lemma \ref{nested} to $T[C_1]$ and $u$ to obtain a sequence $A_i,B_i$ of length $t(\frac{n}{2})+1$ such that $T[A_i],T[B_i]$ are connected, all $B_i$s contain $u$, and the $B_i$s are nested increasing. Then as all $B_i$s contain $u$, $T[B_i\cup C_2]$ is also connected, and thus $(A_i, B_i\cup C_2)$ is a sequence of distinct-size, connected edge-partitions of length at least $t(\frac{n}{2})+1= 2+t(\lceil\frac{n-1}{2}\rceil)-1\ge t(n)-1$ and in case $n$ is of the form $4\cdot 3^\ell+\sum_{i=0}^\ell 3^i$, then $P(T(\ell),2)=t(n)-2$ by Proposition \ref{t}.

\vskip 0.2truecm

\textsc{Case II.} For every edge $uv$ of $T$ one of the components of $T-uv$ is larger than the other.

\vskip 0.15truecm

Orient every edge $uv$ of $T$ towards the larger component of $T-uv$. Clearly, this oriented tree will have exactly one sink $v$ and all components $C_1,C_2,\dots,C_d$ of $T - v$ have size at most $\frac{n-1}{2}$. Suppose the $C_i$s are in decreasing order of size. If $C_1$ has size at least $\frac{n-1}{3}$, then let $T'=T[C_1\cup v]$. If $|C_1|<\frac{n-1}{3}$, then consider the smallest $j$ such that $S:=\sum_{i=1}^j|C_i|\ge \frac{n-1}{3}$. As all $C_i$s are smaller than $\frac{n-1}{3}$, we have $\frac{n-1}{3}\le S\le 2\frac{n-1}{3}$. Therefore either $T[v \cup \bigcup_{i=1}^jC_i]$ or $T[v \cup \bigcup_{i=j+1}^dC_i]$ has its number of vertices between $\lceil \frac{n-1}{3}\rceil$ and $\frac{n}{2}$. We let $T'$ be the one with this property. We apply Lemma \ref{nested} to $T'$ and $v$, and extend the obtained $B_i$ with $V(T)\setminus V(T')$ to reach at least $t(\lceil\frac{n-1}{3}\rceil)+1\ge t(n)-2$ distinct-size connected edge-partitions which is at least $P(T(\ell),2)$ if $n$ is of the form $4\cdot 3^\ell+\sum_{i=0}^\ell 3^i$.
\end{proof}

\begin{proof}[\textbf{Proof of Theorem \ref{constr}}]
The same construction works for all $k\ge 2$. Our construction $G$ contains a complete binary tree $T$ of height $h_1+h_2$ as a spanning tree, so the total number of vertices is $2^{h_1+h_2+1}-1$. Then as a final step, for any vertex $v\in T$ at distance $h_1$ from the root of $T$, we turn the subtree $T_v$ of all descendants of $v$ into a complete subgraph.
In this step, $\binom{2^{h_2+1}-1}{2} - (2^{h_2+1}-2)$ additional edges are inserted inside each of the $2^{h_1}$ subtrees.
The total number of edges is $2^{h_1+1}-2+2^{h_1}\cdot\binom{2^{h_2+1}-1}{2}$, which, if $h_2$ tends to infinity, is $(1+o(1))\cdot n\cdot 2^{h_2}$, and thus the average degree $d$ is approximately $2\cdot 2^{h_2}$.

We prove that $P(G,k)$ is at most $2^{k^2}\cdot\sum_{i=0}^{k-1}2^{i\cdot (2h_2+1)}h_1^{k-1-i}$.
Consider a connected $k$-edge partition $E_1,E_2,\dots,E_k$ of $G$ and let $v_i$ be the unique vertex of $E_i$ that is closest to the root of $T$ in $T$. We assume that the distance of $v_1,v_2,\dots,v_k$ to the root is monotone decreasing. Let $i$ denote the number of $E_j$s with $v_j$ being at a distance at least $h_1$ from the root, i.e. that is completely contained in one of the complete subgraphs at the ``bottom part'' of $G$. As there is at least one part $E_j$ that contains other vertices, we have $i\le k-1$. We have (at most) $2^{i\cdot (2h_2+1)}$ ways to choose the sizes of ``bottom parts'' $E_j$. All ``non-bottom'' part sizes $E_j$ can be determined if all sizes $E_i$ $i<j$ are determined by the following choices: 
the distance between the root and the vertex $v_j$ ($h_1$ choices); whether both children of $w_j$ are incident to $E_j$ or just one ($2$ choices); and what are the part sizes ``below'' $E_j$ that is for which $v_i \in T_{v_j}$ ($2^{j-1}\le 2^{k-1}$ choices). This is a total of at most $2\cdot 2^{k-1}h_1$ possibilities. As once the sizes of $E_1,\dots,E_{k-1}$ are determined, then so is the size of $E_k$, therefore $P(G,k)\le 2^{k^2}\sum_{i=0}^{k-1}2^{i\cdot (2h_2+1)}h_1^{k-1-i}$ as claimed. If $h_1\ge 2^{2h_2}=\Theta(d^2)$, then the summands are maximized when $i=0$ and so $P(G,k)$ is of order $h_1^{k-1}=\Theta(\log^{k-1}n)$, while if $h_1=o(2^{2h_2})=o(d^2)$, then the largest summand corresponds to $i=k-1$ and so $P(G,k)=\Theta(d^{2(k-1)})$. Furthermore, one can set $h_2:=\log_2 d(a)$ and $h_1:=\log_2 n_a-\log_2 d(a)$ to obtain that the average degree of $G$ is $\Theta(d(a))$, the number of vertices is $n_a$, and $P(G,k)=\Theta(d(a)^{2(k-1)})$ if $d(a)=\Omega(\log n_a)$ and $P(G,k)=\Theta(\log_2^{k-1}n)$ if $d(a)=o(\log_2 n_a)$. In the most interesting case, when $d(a)=\omega(\log n_a)$ holds, then we have $P(G,k)=(1+o(1))(d^2)^{k-1}$. This completes the proof.
\end{proof}

\begin{proof}[\textbf{First proof of Theorem \ref{lowerpg}}] Let $t \ge k$ be the  integer such that $8(t+1) > d = d(G) \ge 8t$. Observe that $t> d/8-1$. 

Then by Theorem \ref{mader}, $G$ contains a $(2t+1)$-connected subgraph $H$ with $d(H) > d(G) - 4t \ge 4t$. Observe that $n(H)-1 \ge d(H)> 4t$, hence $n(H) \ge 4t +2 > 4(d/8 -1) +2 = d/2-2$.

Then by Theorem \ref{tutte}, $H$ contains (more than) $t$ edge-disjoint spanning trees, hence we can choose $k$ of them. Let these $k$ edge-disjoint spanning trees be $T_i$ ($i=1,2,\ldots,k$). The number of edges in $H \setminus (T_1 \cup \cdots \cup T_k)$ is $$e(H)-k(n(H)-1)=\frac{d(H)n(H) - 2kn(H) + 2k}{2} > \frac{4tn(H) - 2kn(H)}{2} = n(H)(2t-k) >$$
$$ > (d/2-2)(d/4-2-k)=d^2/8-d(k+3)/2 + 2k + 4.$$


\vspace{2mm}

Now let us define the $k$ parts in the following way: partition the edges in $E(H) \setminus (E(T_1) \cup \cdots \cup E(T_k))$ into parts $A_1, A_2, \ldots , A_k$ with $|A_1| \le |A_2| \le \ldots \le |A_k|$. Then let $E_i := E(T_i) \cup A_i$ for $i=1,2, \ldots, k-1$ and $E_k:= E(T_k) \cup A_k \cup (E(G) \setminus E(H))$. Note that all parts are connected as the $T_i$s are spanning trees in $H$.

Now consider all components of $G \setminus H$. Since $G$ is connected, any component of $G \setminus H$ has a vertex in it adjacent to some vertex in $H$ and hence to a vertex in $T_k$ which is a spanning tree of $H$.

Hence we can add all the components of $G \setminus H$ to $E_k$. This addition to $E_k$ is constant: it depends only on $T_k$ which is fixed but does not depend on $A_k$, therefore the original partitions are kept distinct.

We are done by plugging into the theorem of Erd\H os and Lehner (\cite{EL}) the number of edges in $H \setminus (T_1 \cup \ldots \cup T_k)$ which are free for partition into $k$ parts. The number of edges in $H \setminus (T_1 \cup \ldots \cup T_k)$ is at least $\frac{d^2}{8}-\frac{d(k+3)}{2} + 2k + 4$ to get $P(G,k) \ge \pi(\frac{d^2}{8}-\frac{d(k+3)}{2} + 2k + 4, k) \ge (1-o(1))\frac{d^{2(k-1)}}{8^{k-1}(k-1)!k!}$.
\end{proof}

Now we give a second proof of Theorem \ref{lowerpg}. It only works in the special case $k=2$, and it gives a weaker constant than the first proof, but it is completely elementary as it does not rely on Theorem \ref{mader}, nor on Theorem \ref{tutte}.

\begin{proof}[\textbf{Second proof of Theorem \ref{lowerpg}} for the case\/ $k=2$]
Let $G=(V,E)$ be a connected graph on $n$ vertices with $\frac{dn}{2}$ edges. It is well-known that $G$ contains a subgraph $G'$ with minimum degree at least $\frac{d}{2}$. Therefore $G'$, and thus $G$, contains a path $P$ with at least $\frac{d}{4}+1$ vertices. Let $u_1,u_2,\dots,u_t$ be the vertices of $P$ with $u_iu_{i+1}$ an edge for all $i$, and $t=\lfloor \frac{d}{4}\rfloor+1$. 

Then, as the minimum degree of $G'$ is at least $\frac{d}{2}$, we have $m:=e(G[V(P),V\setminus V(P)])\ge (t+1)(d/2-t)\ge \frac{d^2}{16}$. Let $e_1,e_2,\dots,e_m$ be an enumeration of $E(G[V(P),V\setminus V(P)])$ such that every $u_iv$ comes before every $u_jv'$ for $i<j$. If $e_\ell=u_iv$, then we define $r(\ell)=i$ and let $P_\ell$ be the subpath of $P$ from $u_1$ to $u_{r(\ell)}$. We are ready to define our partitions: let $E^1_\ell$ consist of $e_1,e_2,\dots,e_\ell$, the edges of $E(G[P_\ell])$, and the edges of those components $C$ of $G\setminus P$ that are adjacent to at least one edge $e_i$ with $i\le \ell$. As $E^1_\ell$ contains the edges of $P_\ell$, $E^1_\ell$ induces a connected subgraph of $G$. The same holds for $E^2_\ell:=E\setminus E^1_\ell$. Indeed, it contains the subpath of $P$ from $u_{r(\ell)}$ to $u_k$ and since $G$ is connected, every component of $G\setminus P$ of which the edges belong to $E^2_\ell$ must have a vertex adjacent to some $e_i$ with $i>\ell$. Clearly, $E^1_\ell\subsetneq E^1_{\ell+1}$ for all $\ell$ and thus we defined at least $m/2\ge \frac{d^2}{32}$ connected partitions corresponding to distinct pairs $(m_1,m_2)$.
\end{proof}

Before giving the proof of Theorem \ref{cmcg}, let us mention the obvious relation $P(G,2)\ge \mathrm{CMC}(G)/2$. If $A,B$ is a connected cut with $e(A,B)=\mathrm{CMC}(G)$, then edges of the cut can be added to both $G[A]$ and to $G[B]$ and so
we can have a connected edge-partition with sizes $(e(G[A])+i,e(G[B])+e(A, B)-i)$ for any $0\le i\le \mathrm{CMC}(G)$ and at most two such partitions may yield the same unordered pair of sizes.

\begin{proof}[\textbf{Proof of Theorem \ref{cmcg}}] First we prove the theorem in the case of $r=2$.

\smallskip 

If $d(G)$ is at most 8, then all we need to prove is that $CDC_r(G)\neq 0$, which is clear as $G$ is connected.
So we can suppose $d(G)/4 \ge 2$. Then by Theorem \ref{mader}, $G$ contains a $d(G)/4$-connected subgraph $H$ with $n(H) > \delta(H) \ge d(G)/4$, where $\delta(H)$ denotes the minimum degree of $H$. Applying Theorem \ref{GyLvertex} (as in particular $H$ is $2$-connected), we partition $V(H)=V(A) \cup V(B)$, $|V(A)|=d(G)/8$, $|V(B)|=n(H)-d(G)/8\ge d(G)/8$, where both $A$ and $B$ induce connected parts.

Consider the number $e(A, B)$ of edges in the cut ($A, B$). Every vertex in $A$ can have at most $d(G)/4$ edges inside $A$, hence must have at least $d(G)/4$ edges connected to vertices in $B$. This means $e(A,B)  \ge d^2(G)/64$.

Since $G$ is connected, any component of $G \setminus H$ is adjacent by at least one edge to $H$, hence by at least one edge to either $A$ or $B$. So we can add it to one of the parts, keeping the connectivity of the parts and without decreasing $e(A, B)$. This completes the proof for $r=2$.


\medskip

To prove the statement in the case of $r \ge 3$ we apply a slight modification of the argument. 

\smallskip 

As $r\ge 3$ is fixed, just as in the case $r=2$, we can assume $d(G)/4\ge r$. Then by Theorem \ref{mader}, $G$ contains $d(G)/4$-connected subgraph $H$ with $n(H) > \delta(H) \ge d(G)/4$.

Applying Theorem \ref{GyLvertex} we partition $V(H)$ into $r$ parts (as in particular $H$ is $r$-connected) with cardinalities $n_1=n_2=\ldots=n_{r-1}=\frac{d(G)}{4r}$ and $n_r=n(H)-(r-1)\frac{d(G)}{4r} \ge \frac{d(G)}{4r}$ to get the connected parts $H_1, \ldots, H_{r-1}, H_r$.

Every vertex in $H_j$ ($j=1,\ldots,r-1$) can have at most $\frac{d(G)}{4r}$ edges inside $H_j$ hence at least $(r-1)\frac{d(G)}{4r}$ edges to the other parts.

The number of edges emanating from the parts $H_j$ ($j=1,\ldots,r-1$) is at least $$\frac{(r-1)\cdot\frac{d(G)}{4r}\cdot \frac{(r-1)d(G)}{4r}}{2} \ge \frac{d^2}{72},$$ as $r \ge 3$.

Since $G$ is connected, any component of $G \setminus H$ is adjacent by at least one edge to $H$, hence by at least one edge to one of $H_1, H_2,\dots H_r$. So we can add it to one of the parts, keeping the connectivity of the parts and without decreasing the number of edges between the $r$ parts. This completes the proof.
\end{proof}

\begin{proof}[\textbf{Proof of Theorem \ref{thmonpi}}] 
Let $G$ be a connected graph on $n$ vertices with average degree $d$. We apply Theorem \ref{mader} to obtain a $(\frac{d}{4}+1)$-connected subgraph $H$ of $G$. In particular, the minimum degree of $H$ is at least $\frac{d}{4}+1$, so $H$ contains a path $P$ on $\frac{d}{4}+2$ vertices. We partition $P$ into $k-1$ parts to obtain vertex-disjoint paths $P_1,P_2,\dots,P_{k-1}$
each containing at least $\lfloor\frac{d}{4(k-1)}\rfloor$ vertices. 

\begin{clm}\label{ordered}
Let \/ $(a_1,a_2,\dots,a_k)$ be an ordered\/ $k$-tuple with\/ $\sum_{i=1}^ka_i=n$, $1\le a_i \le \lfloor\frac{d}{4(k-1)}\rfloor$ for all $1\le i\le k-1$. Then there exists a connected\/ $k$-partition\/ $A_1,A_2,\dots,A_k$ of\/ $V(G)$ such that\/ $|A_i\cap P_i|=a_i$ for all\/ $1\le i\le k-1$ and the\/ $k$-tuples\/ $(|A_1|,|A_2|,\dots,|A_k|)$ are pairwise distinct.
\end{clm}

\begin{proof}[Proof of Claim]
Let $(a_1,a_2,\dots,a_k)$ be a $k$-tuple as in the statement of the claim. For $i=1,2,\dots,k-1$, let $X_i$ be the set of the first $a_i$ vertices of $P_i$ and let $X_k=V(H)\setminus \cup_{i=1}^{k-1}X_i$. Observe that $|\cup_{i=1}^{k-1}X_i|\le \frac{d}{4}$, therefore by the connectivity of $H$, we have that $H[X_k]=G[X_k]$ is connected. For any component $C$ of $G-H$, let $r(C)=\min \{j: G[X_j\cup C]~\text{is connected}\}$. As $G$ is connected, $r(C)$ is well-defined for all components. Finally, we define $A_i=X_i\cup \bigcup_{C:r(C)=i}C$ for all $i=1,2,\dots,k$.

By definition, we have $|A_i\cap P_i|=|X_i|=a_i$. Consider two $k$-tuples $(a_1,a_2,\dots,a_k)$, $(a'_1,a'_2,\dots, a'_k)$ and the corresponding partitions $(A_1,A_2,\dots,A_k)$, $(A'_1,A'_2,\dots,A'_k)$. Let $j=\min\{j:a_j\neq a'_j\}$. By symmetry, we can assume $a_j<a'_j$. Then for any $i<j$, we have $X_i=X'_i$, and therefore $A_i=A'_i$. On the other hand, $X_j\subsetneq X'_j$ and as $A_i=A'_i$ for all previous $i$, we have $A_j\subseteq A'_j$. So $(|A_1|,|A_2|,\dots,|A_k|)\neq (|A'_1|,|A'_2|,\dots,|A'_k|)$, proving the claim.
\end{proof}

Every unordered $k$-tuple of sizes of a connected partition can belong to at most $k!$ ordered $k$-tuples, therefore by Claim\ref{ordered}, we have $\pi(G,k)\ge \frac{1}{k!}\cdot \lfloor\frac{d}{4(k-1)}\rfloor^{k-1}$. This completes the proof of the theorem.
\end{proof}

\section{Problems and remarks}

1/ Mader's theorem: We heavily used the theorem of Mader \cite{M72} from 1972 concerning the existence of $k$-vertex-connected subgraphs. However there were succesive improvements by Yuster \cite{Y03}, Bernshteyn and Kostochka \cite{BK16}, Xu, Lai and Tian \cite{XLT20} until Carmesin \cite{C20} proved the best possible bound: every graph with $d(G) \ge (3+ \frac{1}{3})k$ has a ($k+1$)-connected subgraph with more than $2k$ vertices.

We have decided to use the classical, easy-to-work-with version of the theorem at the cost of getting slightly weaker bounds up to a small multiplicative factor.
We also mention here the $k$-edge-connected theorem of Mader \cite{M74}, see also \cite{GGL95}, for the sake of completeness and as a tool in an alternative approach: Let $G$ be a graph on $n$ vertices and with $e(G) \ge kn - \binom{k+1}{2} + 1$. Then $G$ has a ($k+1$)-edge-connected subgraph, and this is the best possible.

\medskip

2/ We proved in Theorem \ref{lowerpg} a lower bound on $P(G,k)$ of order $d(G)^{2(k-1)}$. Essentially the same proof can be adopted when $d(G)$ is replaced by the maximum-average degree $Mad(G):=\max\{d(H): H \textrm{ is induced subgraph of }G\}$.

This is because $Mad(G) \ge d(G)$, hence in the theorem of Mader if $Mad(G) \ge 4k$, then we can apply it directly on the induced subgraph $H$ of $G$, where $Mad(G)=d(H)$ is realized to obtain the ($k+1$)-connected subgraph. 

\medskip

3/ We proved in Theorem \ref{cmcg} a lower bound $\mathrm{CMC}_r(G) = \Omega_r (d^2)$ on the connected color classes $r$-partite problem (motivated by the Max Cut and the maximum $r$-partite problem).

This seems as far as we are aware the first time such bounds are given and up to constant multiplicative factors are sharp. In this case, we can also develop a parallel argument replacing $d(G)$ by $Mad(G)$.
It would be interesting to obtain for $r \ge 2$, the best possible constants $c_r$ in Theorem \ref{cmcg}.

\medskip

4/ One might wonder whether we need the extra constant multiplicative factor for general connected graphs compared to the complete ternary tree $T(n,3)$ in Theorem \ref{kmain}. 

\begin{conjecture}
$P(G,2)$ is  minimized by\/ $T(n,3)$ over all connected graphs on\/ $n$ vertices whenever\/ $n=\sum_{i=0}^\ell3^i$ for some\/ $\ell$.
\end{conjecture} 

The proof might be hard as $P(G,2)$ is not monotone in $G$. 

\begin{example}
Let\/ $G$ be the graph that we obtain from the complete binary tree on\/ $31$ vertices by adding\/ $8$ edges joining neighboring leaves. Formally, let\/ $V(G)=\{1,2,\dots,31\}$ and\/ $E(G)=\{(i,2i), (i,2i+1): i=1,2,\dots,15\}\cup \{(i,i+1): i=16,18,20,22,24,26,28,30\}$; and let\/ $G'=G\setminus e$ with\/ $e=(30,31)$. Then\/ $P(G,2)=8$ as there exist connected edge sets of size\/ $1,2,3,4,8,9,18$ and\/ $19$, while\/ $P(G\setminus e,2)=9$ as there we can have connected edge sets of size\/ $1,2,3,4,7,8,9,17$, and\/ $18$. 
\end{example}

\noindent \textbf{Acknowledgement}. We would like to thank the unknown referees for their careful reading and their many suggestions that helped us improve the presentation of this paper.

\end{document}